\documentclass[11pt]{article}
\usepackage{amssymb}
\usepackage{latexsym}
\usepackage[dvips]{epsfig}

\newcommand{\qed}{\mbox{$\Diamond$}\vspace{\baselineskip}}

\newtheorem{theorem}{Theorem}

\newcommand{\cN}{{\cal N}}

\newcommand{\la}{\lambda}

\newcommand{\vanish}[1]{}
\begin{document}
  
\setlength{\baselineskip}{1.2\baselineskip}
\title{Two Injective Proofs of a Conjecture of Simion}
\author{Mikl\'os B\'ona \\
        Department of Mathematics \\
        University of Florida\\
        Gainesville, FL 32611\\
	USA\\
	bona@math.ufl.edu\\[5pt]
	and\\[5pt]
	Bruce E. Sagan\\
	Department of Mathematics\\ 
	Michigan State University\\
	East Lansing, MI 48824-1027\\
	USA\\
	sagan@math.msu.edu}
\maketitle

\begin{abstract}
Simion~\cite{simion} conjectured the unimodality of a sequence counting
lattice paths 
in a grid with a Ferrers diagram removed from the northwest corner.
Recently, Hildebrand~\cite{hildebrand} and then Wang~\cite{wan:spc}
proved the stronger result that this sequence is actually log concave.
Both proofs were mainly algebraic in nature.
We give two combinatorial proofs of this theorem.
\end{abstract}

\section{Introduction}
In this note we present two injective proofs of a strengthening of a
conjecture of Simion~\cite{simion}.  To describe the result,
let $\lambda=(\lambda_1, \lambda_2, \cdots ,\lambda_k)$ be the Ferrers
diagram of a partition viewed as a set of squares in English notation.
(See any of the texts~\cite{bona,sag:tsg,sta:ec2} for definitions of
terms that we do not define here.)  The shape $\lambda$ will be
fixed for the rest of this paper.

\begin{figure}[ht]
 \begin{center}
  \epsfig{file=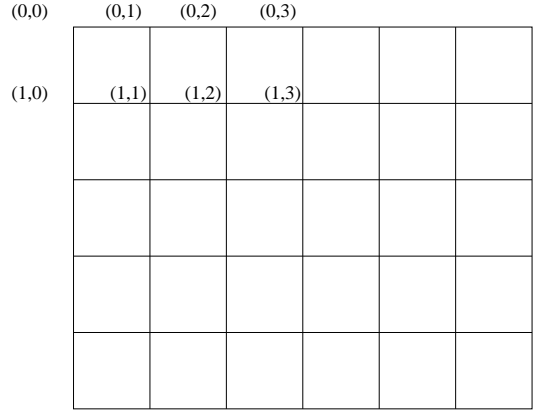}
  \caption{Labeling our points in the grid. }
  \label{rgrid}
 \end{center}
\end{figure}

Consider a grid with the vertices labeled $(i,j)$ for $i,j\ge0$ as in
Figure~\ref{rgrid}.
Place $\lambda$ in the northwest corner of this array so that its
squares coincide with those of the grid.

A {\em northeastern lattice path\/} is a lattice path on the grid in
which each step goes one unit to the north or one unit to the east. 
Let $N(m,n)$ be the number of northeastern lattice paths from $(m,0)$ to 
$(0,n)$ that {\em do not go inside $\lambda$\/} (although they may
touch its boundary), and let $\cN(m,n)$ be the set of such paths.
In particular, $N(m,n)=0$ if either the starting or ending point is
inside $\lambda$.

Simion~\cite{simion} conjectured that for all $m,n\ge0$ the sequence
$$
N(0,m+n), N(1,m+n-1), \ldots, N(m+n,0)
$$
is unimodal.  Lattice path techniques  for proving unimodality were
investigated by Sagan~\cite{sagan}, but the conjecture remained open
at that point.  Recently, Hildebrand~\cite{hildebrand} proved the
stronger result that this sequence is actually log concave by mostly
algebraic means.  Shortly thereafter, Wang~\cite{wan:spc} simplified
Hildebrand's proof using results about Polya frequency sequences.  In
the present work, we will give two injective proofs of the stong
version of Simion's conjecture.  The one in Sections~\ref{seceasy}
and~\ref{sechard} employs ideas from
Hildebrand's proof while the one in Section~\ref{direct} is more
direct.  Our injections come 
from a method of Lindstr\"om~\cite{lin:vri}, later
popularized by Gessel and Viennot~\cite{gv:bdp,gv:dpp}, that can be
used to prove total positivity results for matrices.  For an
exposition, see Sagan's book~\cite[pp.\ 158--163]{sag:tsg}.
B\'ona~\cite{logc} has used related ideas to prove the log concavity
of a sequence counting $t$-stack sortable permutations.

We end this section by reiterating the statement of the main theorem for
easy reference.  Notice that when $\lambda=\emptyset$ it specializes to
the well-known result that the rows of Pascal's triangle are log concave.
\begin{theorem}[The Strong Simion Conjecture]
\label{main}
Let $\lambda$ be the Ferrers diagram of a partition and let $N(m,n)$ be
the number of northeastern lattice paths 
in the grid from $(m,0)$ to $(0,n)$  which do not
intersect the interior of $\lambda$.  Then for all $m,n\ge0$ the sequence
$$
N(0,m+n), N(1,m+n-1), \ldots, N(m+n,0)
$$
is log concave.
\end{theorem}

\section{A decomposition of the problem}

This preliminary part of the first proof is from~\cite{hildebrand}.
We include it so that our exposition will be self
contained.  We need to prove that for all $m,n>0$ we have
$$
N(m-1,n+1) N(m+1,n-1) \leq N(m,n)^2.
$$
To prove this, it suffices to show that 
$$
N(m-1,n+1)N(m+1,n) \leq N(m,n) N(m,n+1),
$$
because then, by symmetry, we also have
$$
N(m+1,n-1)N(m,n+1) \leq N(m,n) N(m+1,n).
$$
Now multiplying the last two equations together and simplifying gives
the first.

The second inequality can be proved by demonstrating another pair of
equations, namely
\begin{equation} 
\label{easy} 
N(m,n+1) N(m+1,n) \leq N(m,n) N(m+1,n+1),
\end{equation} 
and
\begin{equation} 
\label{hard}
 N(m-1,n+1)N(m+1,n+1)   \leq  N(m,n+1)^2.
\end{equation}
Multiplying these two equations together and cancelling gives the
desired result.

\section{The proof of (\ref{easy})}		
\label{seceasy}

In this section we prove  that (\ref{easy}) holds by constructing an
injection 
$$
\Psi: \cN (m,n+1) \times \cN (m+1,n) \rightarrow\cN(m,n) \times \cN(m+1,n+1).
$$

\begin{figure}[ht]
 \begin{center}
  \epsfig{file=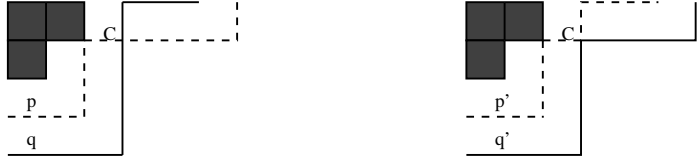}
  \caption{The action of $\Psi$. }
  \label{easyf}
 \end{center}
\end{figure}

Consider a path pair $(p,q)\in \cN (m,n+1) \times \cN (m+1,n)$. Then $p$
and $q$ must intersect. Let $C$ be their first (most southwestern)
intersection point. Say that $C$ splits $p$ into parts $p_1$ and $p_2$,
and splits $q$ into parts $q_1$ and $q_2$. Then the
concatenation of $p_1$ and $q_2$ is a path in $\cN(m,n)$, and the
concatenation of $q_1$ and $p_2$ is a path in $\cN(m+1,n+1)$. So
define $\Psi(p,q)=(p_1q_2,q_1p_2)=(p',q')$.  It is easy to see that the image
of $\Psi$ is exactly all 
$(p',q')\in \cN(m,n) \times \cN(m+1,n+1)$ such that $p'$ and $q'$
intersect.  It is also simple to verify that if $\Psi(p,q)=(p',q')$, 
then applying the same algorithm to $(p',q')$ recovers $(p,q)$.  So
$\Psi$ is injective.
See Figure \ref{easyf} for an example. 

\section{The proof of (\ref{hard})}
\label{sechard}

In this section we  construct an injection
$$
\Phi:  \cN(m-1,n)\times\cN(m+1,n)   \rightarrow  \cN(m,n)^2,
$$
thus proving~(\ref{hard}).

Let $(p,q)\in  \cN (m-1,n) \times \cN (m+1,n)$.  If $P=(i,j)$ and
$Q=(i,k)$ are vertices of $p$ and $q$, respectively, with the same first
coordinate, define the {\it vertical distance from $p$ to $q$ at
$P$ and $Q$\/} to be $k-j$.  The vertical distance from $p$
to $q$ starts at 2 for their initial vertices and ends at 0 for their
final ones.  Since vertical distance can change by at most one with a
step of a path, there must be some vertical distance equal to 1.  Let
$P$ and $Q$ be the first (most southwest) pair of points with vertical
distance one.  

\begin{figure}[ht]
 \begin{center}
  \epsfig{file=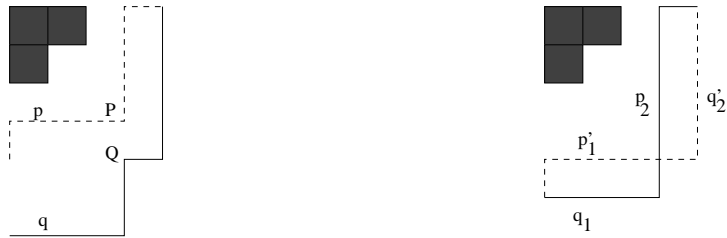}
  \caption{The action of $\Phi$. }
  \label{newphi}
 \end{center}
\end{figure}

It follows from our choice of vertices that $p$ must
enter $P$ with an east step and $q$ must enter $Q$ with a north step.  Let
$p_1$ and $p_2$ be the portions of $p$ before and after $P$,
respectively, and similarly for $q$.  Now let
$$
\Phi(p,q)=(p_1'q_2,q_1'p_2)
$$ 
where $p_1'$ is $p_1$ moved south one
unit and $q_1'$ is $q_1$ moved north one unit.  Since $P$ and $Q$ are
the first pair of points at vertical distance one, $q_1'$ will not
intersect $\la$ and the concatenations are valid paths in $\cN(m,,n)$.
In fact, the image of $\Phi$ is exactly all path pairs
$(p',q')\in\cN(m,n)^2$ that have a pair of points at vertical distance
-1.  Applying the same procedure to the first such pair inverts the
map and so $\Phi$ is injective.  See Figure \ref{newphi} for an example.

This completes the first proof of Theorem~\ref{main}. \quad\qed

\section{A more direct proof}
\label{direct}

The reader may wonder if we can do away with splitting our problem into
two parts, that is, equations (\ref{easy}) and (\ref{hard}). The answer is
yes, and the necessary injection is just a modification $\overline{\Phi}$ of the map
$\Phi$.  This will give us a second, completely combinatorial, proof
of our main theorem

Take a path pair $(p,q)\in \cN(m-1,n+1)\times\cN(m+1,n-1)$.  Notice
that $p$ and $q$ must intersect.  So before the first intersection
there must be a first pair of points $P, Q$ (on $p,q$ respectively) at
vertical distance 1.  Similarly, after the last intersection there
must be a last pair of points $\overline{P},\overline{Q}$ at horizontal distance
1, where horizontal distance is defined analogously.  Let $P$ and
$\overline{P}$ divide $p$ into subpaths $p_1,p_2,p_3$ and use the same
notation for $q$.  Then define
$$
\overline{\Phi}(p,q)=(p_1'q_2,p_3'',q_1'p_2q_3'')
$$
where $p_1'$ is $p_1$ moved south one unit, $p_3''$ is $p_3$ moved
west one unit, and $q_1',q_3''$ are defined in the analogous way but
moving in the opposite directions.  It is a simple job to verify that
$\overline{\Phi}$ is well-defined and injective just as we did with $\Phi$.

This completes the second proof of Theorem~\ref{main}.\quad\qed

We have two final remarks.  First of all, it is clear from the
geometry of the situation that if $\la$ is self-conjugate then the
sequence in Theorem~\ref{main} is also symmetric, but this does not
hold in general.  One might also wonder if this sequence has the
stronger property that the associated polynomial generating function
has only real zeros.  This is not always true as can be seen by taking
$\la=(1)$ and $m+n=4$.  In this case the associated polynomial is
$x(3x^2+5x+3)$ which has two complex roots.  It might be interesting
to determine for which shapes the real zero property holds.

\end{document}